\documentclass[12pt]{article}
\usepackage{amssymb,amsthm,amsmath,latexsym}
\usepackage{url}
\newtheorem{thm}{Theorem}
\newtheorem{prop}[thm]{Proposition}

\theoremstyle{remark}
\newtheorem{rem}[thm]{Remark}
\newcommand{\FF}{\mathbb{F}}
\newcommand{\ZZ}{\mathbb{Z}}

\DeclareMathOperator{\wt}{wt}

\DeclareMathOperator{\supp}{supp}

\begin{document}
\title{New doubly even self-dual codes having minimum weight 20}

\author{
Masaaki Harada\thanks{
Research Center for Pure and Applied Mathematics,
Graduate School of Information Sciences,
Tohoku University, Sendai 980--8579, Japan.
email: {\tt mharada@tohoku.ac.jp}.}
}
\date{}

\maketitle

\begin{abstract}
In this note, we construct new doubly even self-dual codes having minimum
weight $20$ for lengths $112$, $120$ and $128$.
This implies that there are at least three inequivalent
extremal doubly even self-dual codes of length $112$.
\end{abstract}

%%%%%%%%%%%%%%%%%%%%%%%%%%%%%%
\section{Introduction}

Self-dual codes are an important class of linear codes for both
theoretical and practical reasons (see~\cite{RS-Handbook}).
It is a fundamental problem to determine the largest minimum weights 
among self-dual codes of that length and to construct self-dual
codes with the largest minimum weight.

Let $\FF_2$ denote the finite field of order $2$.
Codes over $\FF_2$ are called {\em binary} and
%% all codes in this paper are binary unless otherwise noted.
all codes in this note are binary.
The \textit{dual code} $C^{\perp}$ of a code 
$C$ of length $n$ is defined as
$
C^{\perp}=
\{x \in \FF_2^n \mid x \cdot y = 0 \text{ for all } y \in C\},
$
where $x \cdot y$ is the standard inner product.
A code $C$ is called 
% \textit{self-orthogonal} if $C \subset C^{\perp}$, 
% and $C$ is  called 
\textit{self-dual} if $C = C^{\perp}$. 
A self-dual code $C$ is {\em doubly even} if all
codewords of $C$ have weight divisible by four, and {\em
singly even} if there is at least one codeword of weight $\equiv 2
\pmod 4$.
It is  known that a self-dual code of length $n$ exists 
if and only if  $n$ is even, and
a doubly even self-dual code of length $n$
exists if and only if $n$ is divisible by eight.
The minimum weight $d$ of a doubly even self-dual code of length $n$
is bounded by
\begin{equation}\label{eq:bound}
d  \le 4  \left\lfloor{\frac {n}{24}} \right\rfloor + 4,
\end{equation}
\cite{MS73}.
A doubly even self-dual code meeting the bound is called  {\em extremal}.
% The minimum weight $d$ of a singly even self-dual code of length $n$
% is bounded by
% $d  \le 4 \lfloor{\frac {n}{24}} \rfloor + 6$ 
% if $n \equiv 22 \pmod {24}$, 
% $d  \le 4  \lfloor{\frac {n}{24}} \rfloor + 4$ 
% otherwise~\cite{Rains}.
% A self-dual code meeting the bound is called  {\em extremal}.

In this note, we study the existence of doubly even self-dual
codes having minimum weight $20$.
By~\eqref{eq:bound}, if there is a doubly even
self-dual code of length $n$ and minimum weight $20$, then
$n \ge 96$.
For length $96$, it is unknown whether there is an
extremal doubly even self-dual code.
For length $104$, the extended quadratic residue code
is the only known extremal doubly even self-dual code 
(see~\cite{RS-Handbook}).
For length $112$, the first extremal doubly even self-dual
code was found in~\cite{H112}.
For lengths $120$ and $128$, it is unknown whether there is an
extremal doubly even self-dual code.
The first doubly even self-dual code 
of length $120$ and minimum weight $20$
was found in~\cite{GNW}.
Then $25$ more doubly even self-dual codes
of length $120$ and minimum weight $20$
were found in~\cite{YW}.
The existence of a doubly even self-dual code
of length $128$ and minimum weight $20$
is known~\cite{G}.
In this note, we construct new doubly even self-dual codes having minimum
weight $20$ for lengths $112$, $120$ and $128$.
This implies that there are at least three inequivalent
extremal doubly even self-dual codes of length $112$.

All computer calculations in this note were
done with the help of {\sc Magma}~\cite{Magma}.

%%%%%%%%%%%%%%%%%%%%%%%%%%%%%%%%%%%%%%%%%%%%%%
\section{Preliminaries}\label{Sec:2}

Let $C$ be a code of length $n$.
The elements of $C$ are called {\em codewords} and the {\em weight}
$\wt(x)$ of a codeword $x$ is the number of non-zero coordinates.
The {\em support} of a codeword $x=(x_1,x_2,\ldots,x_n)$ is 
$\{i \mid x_i=1\}$.
We denote the support of $x$ by  $\supp(x)$.
The minimum non-zero weight of all codewords in $C$ is called 
the {\em minimum weight} of $C$.
Let $A_i$ be the number of codewords of
weight $i$ in $C$.
The {\em weight enumerator} $W_C$ of $C$ is given by
$\sum_{i=0}^n A_i y^i$.
%An $[n,k,d]$ code is an $[n,k]$ code with minimum (non-zero) weight $d$.
By Gleason's theorem~\cite{Gleason}  (see also~\cite{MS73}),
the weight enumerator $W_C$ of a self-dual code of length $n$
is written as:
\begin{equation}\label{eq:W}
W_C = \sum_{j=0}^{\lfloor n/8 \rfloor} a_j(1+y^2)^{n/2-4j}(y^2(1-y^2)^2)^j,  
\end{equation}
for some integers $a_j$.
In addition, 
a doubly even self-dual code of length $n$
exists then $n$ is divisible by eight, and 
the weight enumerator $W_C$ of a doubly even self-dual code of length $n$
is written as:
\begin{equation}\label{eq:WII}
W_C = \sum_{j=0}^{\lfloor n/24 \rfloor} 
a_j(1+14y^4+y^8)^{n/8-3j}(y^4(1-y^4)^4)^j, 
\end{equation}
for some integers $a_j$.
Then Mallows and Sloane~\cite{MS73} established the upper bound
\eqref{eq:bound} on the minimum weights of doubly even self-dual codes.

Let $C$ be a singly even self-dual code and
let $C_0$ denote the
subcode of codewords having weight $\equiv0\pmod4$.
Then $C_0$ is a subcode of codimension $1$.
The {\em shadow} $S$ of $C$ is defined to be 
$C_0^\perp \setminus C$~\cite{C-S}.
There are cosets $C_1,C_2,C_3$ of $C_0$ such that
$C_0^\perp = C_0 \cup C_1 \cup C_2 \cup C_3 $, where
$C = C_0  \cup C_2$ and $S = C_1 \cup C_3$.
If $C$ is a singly even
self-dual code of length divisible by $8$, then $C$ has two doubly
even self-dual neighbors, namely $C_0 \cup C_1$ and $C_0 \cup C_3$
(see~\cite{BP}).
Let $B_i$ be the number of vectors of
weight $i$ in $S$.
The weight enumerator $W_S$ of $S$ is given by
$\sum_{i=d(S)}^{n-d(S)} B_i y^i$, respectively,
where $d(S)$ denotes the minimum weight of $S$.
If $W_C$ is written as in \eqref{eq:W}, then $W_S$ can
be written as follows~\cite[Theorem~5]{C-S}:
\begin{equation}\label{eq:WS}
W_S = \sum_{j=0}^{\lfloor n/8 \rfloor}
(-1)^ja_j2^{n/2-6j}y^{n/2-4j}(1-y^4)^{2j}.
\end{equation}

Two self-dual codes $C$ and $C'$ of length $n$
are said to be {\em neighbors} if $\dim(C \cap C')=n/2-1$. 
% Any self-dual code of length $n$ can be reached
% from any other by taking successive neighbors (see~\cite{C-S}).
% It is known that a self-dual code $C$ of length $n$ has 
% $2(2^{n/2-1}-1)$ self-dual neighbors.
% These neighbors are constructed by finding
% $2^{n/2-1}-1$ subcodes of codimension $1$ in $C$
% containing the all-one vector.
%% Two codes $C$ and $C'$ are {\em equivalent}, denoted $C \cong C'$,
Two codes are {\em equivalent} 
if one can be
obtained from the other by permuting the coordinates.
An {\em automorphism} of a code $C$ is a permutation of the coordinates of $C$
which preserves $C$.
The set consisting of all automorphisms of $C$ is called the
{\em automorphism group} of $C$.
% and it is denoted by $\Aut(C)$. 

An $n \times n$ circulant matrix has the following form:
\[
\left(
\begin{array}{ccccc}
r_0&r_1&r_2& \cdots &r_{n-1} \\
r_{n-1}&r_0&r_1& \cdots &r_{n-2} \\
\vdots &\vdots & \vdots && \vdots\\
r_1&r_2&r_3& \cdots&r_0
\end{array}
\right),
\]
so that each successive row is a cyclic shift of the previous one.
Let $A$ and $B$ be $n \times n$ circulant matrices.
Let $C$ be a code with generator matrix of the following form:
\begin{equation} \label{eq:GM}
\left(
\begin{array}{ccc@{}c}
\quad & {\Large I_{2n}} & \quad &
\begin{array}{cc}
A & B \\
B^T & A^T
\end{array}
\end{array}
\right),
\end{equation}
where $I_n$ denotes the identity matrix of order $n$
and $A^T$ denotes the transpose of a matrix $A$.
It is easy to see that $C$ is self-dual if
$AA^T+BB^T=I_n$.
The codes with generator matrices of the form~\eqref{eq:GM}
are called {\em four-circulant}~\cite{4cir}.
In this note, 
we found a 
singly even self-dual four-circulant code 
of length $112$ and minimum weight $18$ and
doubly even self-dual four-circulant codes 
of length $n$ and minimum weight $20$ for $n=112,120,128$,
by a non-exhaustive search.
An exhaustive search is beyond our current computer resources.

%%%%%%%%%%%%%%%%%%%%%%%%%%%%%%%%%%%%%%%%%%%%%%
%%%%%%%%%%%%%%%%%%%%%%%%%%%%%%%%%%%%%%%%%%%%%%
\section{New extremal doubly even self-dual codes of length 112}

\subsection{A singly even self-dual code
of length 112 and minimum weight 18}

By a non-exhaustive search, we found a 
singly even self-dual four-circulant code $C_{112}$
of length $112$ and minimum weight $18$.
The first rows $r_A$ and $r_B$ of $A$ and $B$
in the generator matrix \eqref{eq:GM}
of $C_{112}$ are as follows:
\begin{align*}
r_A&=(1000010101101101111011011010),\\
r_B&=(0010001110000110001010000001), 
\end{align*}
respectively.

Let $C$ be a singly even self-dual
code of length $112$ and minimum weight $18$.
Let $S$ be the shadow of $C$.
% By~\cite[Theorem~5]{C-S},
From \eqref{eq:W} and \eqref{eq:WS},
the possible weight enumerators of $C$ and $S$
are determined as follows:
\begin{align*}
W_{112}^C=& 
1 
+( 99176 + a)y^{18} 
+( 355740 + 16 b + 2 a)y^{20} 
\\&
+( 1745240 + 1024 c - 64 b - 17 a)y^{22} 
\\&
+( 44404374 + 65536 d - 10240 c - 160 b - 36 a)y^{24} 
\\&
+( 572977944 - 4194304 e - 1048576 d + 33792 c + 960 b + 135 a)y^{26}
\\&
+ \cdots,
\\
W_{112}^S=& 
ey^4 
+(- 26 e + d)y^8 
+( 325 e - 24 d - c)y^{12} 
\\&
+(- 2600 e + 276 d + 22 c + b)y^{16} 
\\&
+( 14950 e - 2024 d - 231 c - 20 b - 4 a)y^{20} 
+ \cdots,
\end{align*}
respectively, where $a,b,c,d,e$ are integers.
In order to determine the weight enumerator of $C_{112}$,
we found that
\[
A_{18}=8512,
d(S)=16,
B_{16}=728.
\]
This gives
\[
a=-90664,
b=728,
c=d=e=0.
\]
The weight distribution of $C_{112}$ is listed in Table~\ref{Tab:WD}.
Note that singly even self-dual
codes of length $112$ and minimum weight $18$ with weight
enumerators corresponding to $e=1$ were found in~\cite{H112}.

%%%%%%%%%%%%%%%%%%%%%%%%%%%%%%
\begin{table}[thb]
\caption{Weight distribution of $C_{112}$}
\label{Tab:WD}
\begin{center}
%{\small
{\footnotesize
%{\scriptsize
\begin{tabular}{c|r|c|r}
\noalign{\hrule height0.8pt}
$i$ & \multicolumn{1}{c|}{$A_i$} & 
$i$ & \multicolumn{1}{c}{$A_i$} \\
\hline
$ 0,112$&                 1& $38, 74$ &    31676520067584 \\
$18, 94$&              8512& $40, 72$ &   109690203298312 \\
$20, 92$&            186060& $42, 70$ &   325630986391040 \\
$22, 90$&           3239936& $44, 68$ &   831288282918576 \\
$24, 88$&          47551798& $46, 66$ &  1829637194737408 \\
$26, 86$&         561437184& $48, 64$ &  3479230392288469 \\
$28, 84$&        5424089452& $50, 62$ &  5725819388994432 \\
$30, 82$&       43459872064& $52, 60$ &  8165553897114152 \\
$32, 80$&      291008417322& $54, 58$ & 10099951175046656 \\
$34, 78$&     1639219687168& $56$ & 10841051388476292 \\
$36, 76$&     7813559379696&     &\\
\noalign{\hrule height0.8pt}
\end{tabular}
}
\end{center}
\end{table}
%%%%%%%%%%%%%%%%%%%%%%%%%%%%%%%

% We verified that the code $C_{112}$ has automorphism group
% of order $112$.

%%%%%%%%%%%%%%%%%%%%%%%%%%%%%%%%%%%%%%%%%%%%%%
\subsection{New extremal doubly even self-dual codes of length 112}

% Two self-dual codes $C$ and $C'$ of length $n$
% are said to be {\em neighbors} if $\dim(C \cap C')=n/2-1$. 
% Any self-dual code of length $n$ can be reached
% from any other by taking successive neighbors (see~\cite{C-S}).
% It is known that a self-dual code $C$ of length $n$ has 
% $2(2^{n/2-1}-1)$ self-dual neighbors.
% These neighbors are constructed by finding
% $2^{n/2-1}-1$ subcodes of codimension $1$ in $C$
% containing the all-one vector.
If $C$ is a singly even
self-dual code of length divisible by $8$, then $C$ has two doubly
even self-dual neighbors 
%, namely $C_0 \cup C_1$ and $C_0 \cup C_3$
(see Section \ref{Sec:2}).
We verified that the two doubly even self-dual neighbors of
$C_{112}$ have minimum weights $20$ and $16$.
We denote the extremal doubly even self-dual neighbor by $D_{112}$.
The code $D_{112}$ is constructed as
\[
\langle (C_{112} \cap \langle x \rangle^\perp), x \rangle,
\]
where 
the support $\supp(x)$ of $x$ is
\[
\{1,2,3,9,10,12,13,16,24,61,66,69,82,89,92,96,97,108,111,112\}.
\]
% We verified that the code $D_{112}$ has automorphism group
% of order $112$.

Moreover, by a non-exhaustive search, we found an extremal 
doubly even self-dual four-circulant code $E_{112}$.
The first rows $r_A$ and $r_B$  of $A$ and $B$
in the generator matrix~\eqref{eq:GM}
of $E_{112}$ are as follows:
\begin{align*}
r_A&=(1000000010001000011111001101),\\
r_B&=(0111101000101001110101011110),
\end{align*}
respectively.
% We verified that the code $E_{112}$ has automorphism group
% of order $56$.

We denote the known extremal
doubly even self-dual code in~\cite{H112} by $H_{112}$.
In order to distinguish $H_{112}, D_{112}$ and $E_{112}$,
we consider the following invariant.
Let $C$ be an extremal doubly even self-dual code of length $112$.
Let $M(C)$ be the matrix with rows composed of the codewords of 
weight $20$ in $C$, 
where the $(1,0)$-matrix $M(C)$ is regarded as a matrix over $\ZZ$.
Let $m_{i,j}$ denote the $(i,j)$-entry of the matrix
$M(C)^T M(C)$.
Then define
\[
m(C)=\{m_{i,j} \mid i,j \in \{1,2,\ldots,112\}\}.
\]
In Table~\ref{Tab:mC}, we list all elements of $m(C)$ for
$C=H_{112},D_{112}$ and $E_{112}$.
Table~\ref{Tab:mC} shows that the three codes are inequivalent.

\begin{prop}
There are at least three inequivalent extremal doubly even self-dual
codes of length $112$. 
\end{prop}

\begin{rem}
The code $H_{112}$ has automorphism group of order $112$~\cite{H112}.
We verified that the codes $D_{112}$ and $E_{112}$
have automorphism group of order $112$.
\end{rem}

%%%%%%%%%%%%%%%%%%%%%%%%%%%%%%
\begin{table}[thb]
\caption{$m(H_{112}),m(D_{112})$ and $m(E_{112})$}
\label{Tab:mC}
\begin{center}
%{\small
%{\footnotesize
{\scriptsize
\begin{tabular}{l}
\noalign{\hrule height0.8pt}
\multicolumn{1}{c}{$m(H_{112})$} \\
\hline
10613, 10649, 10661, 10703, 10709, 10715, 10721, 10727, 10733, 10739, 10745, \\
10769, 10775, 10781, 10787, 10799, 10805, 10811, 10823, 10829, 10835, 10841, \\
10847, 10853, 10859, 10865, 10871, 10883, 10895, 10901, 10907, 10913, 10919, \\
10925, 10931, 10937, 10943, 10949, 10967, 10973, 10985, 10991, 10997, 11009, \\
11021, 11033, 11045, 11057, 11063, 11069, 11093, 11099, 11117, 63525\\
\hline\hline
\multicolumn{1}{c}{$m(D_{112})$} \\
\hline
10618, 10663, 10672, 10702, 10708, 10717, 10735, 10750, 10765, 10768, 10771, \\
10777, 10783, 10786, 10789, 10801, 10810, 10819, 10831, 10834, 10837, 10840, \\
10843, 10846, 10849, 10852, 10858, 10861, 10864, 10867, 10873, 10882, 10885, \\
10900, 10903, 10906, 10909, 10912, 10918, 10921, 10924, 10927, 10930, 10936, \\
10945, 10954, 10957, 10978, 10984, 10987, 11002, 11011, 11023, 11041, 11044, \\
11056, 11065, 11080, 11086, 11098, 11110, 63525\\
\hline\hline
\multicolumn{1}{c}{$m(E_{112})$} \\
\hline
10581, 10620, 10641, 10653, 10659, 10668, 10674, 10689, 10698, 10701, 10704, \\
10707, 10719, 10728, 10734, 10749, 10758, 10761, 10764, 10770, 10776, 10779, \\
10782, 10785, 10791, 10794, 10797, 10806, 10809, 10812, 10815, 10818, 10821, \\
10824, 10827, 10830, 10833, 10842, 10848, 10851, 10854, 10860, 10863, 10866, \\
10872, 10875, 10878, 10881, 10884, 10890, 10893, 10896, 10899, 10902, 10905, \\
10911, 10914, 10917, 10923, 10926, 10929, 10932, 10935, 10938, 10941, 10950, \\
10953, 10959, 10965, 10968, 10971, 10977, 10980, 10989, 10992, 10995, 11001, \\
11013, 11016, 11025, 11028, 11040, 11046, 11049, 11052, 11055, 11073, 11085, \\
11103, 11151, 63525\\
\noalign{\hrule height0.8pt}
\end{tabular}
}
\end{center}
\end{table}
%%%%%%%%%%%%%%%%%%%%%%%%%%%%%%%

% The covering radius $R(C)$ of a code $C$ of length $n$ is the 
% smallest integer $R$
% such that spheres of radius $R$ around codewords of $C$ cover the
% space $\FF_2^n$.
% By the Delsarte bound, the covering radius of an 
% extremal doubly even self-dual code of length $112$
% is at most $20$. 
% We may suppose without loss of generality that
% $D_{112}= (C_{112})_0 \cup (C_{112})_1$.
% Since $(C_{112})_2 \cup (C_{112})_3$ is a coset of $D_{112}$
% and it has minimum weight $18$, <--- not 18, it is 16!
% we have that $18 \le R(D_{112}) \le 20$.
% For small lengths, the Delsarte bound
% seems to give a rather good upper bound on the covering
% radii of extremal doubly even self-dual codes.
% However, the covering radii of many
% extremal doubly even self-dual codes of length $64$
% do not meet the Delsarte bound, and $64$ seems to be 
% the smallest length for which the Delsarte bound is
% rarely met~\cite{HM-CR}.

%%%%%%%%%%%%%%%%%%%%%%%%%%%%%%%%%%%%%%%%%%%%%%
\section{New doubly even self-dual codes
of length 120 and minimum weight 20}

% By Gleason's theorem (see~\cite{MS73}),
From \eqref{eq:WII},
the possible weight enumerator of a doubly even 
self-dual code of length $120$ and minimum weight $20$
is determined as follows:
\begin{align*}
&
1 
+ a y^{20} 
+( 39703755 - 20 a )y^{24} 
+( 6101289120 + 190 a )y^{28} 
\\&
+( 475644139425 - 1140 a )y^{32} 
+( 18824510698240  + 4845 a )y^{36} 
\\&
+( 397450513031544 - 15504 a )y^{40} 
+( 4630512364732800 + 38760 a )y^{44} 
\\&
+( 30531599026535880  - 77520 a )y^{48} 
+( 116023977311397120  + 125970 a )y^{52} 
\\&
+( 257257766776517715  - 167960 a )y^{56} 
\\&
+( 335200280030755776  + 184756 a )y^{60}
+ \cdots, 
\end{align*}
where $a$ is the number of codewords of weight $20$.

The first doubly even self-dual code
of length $120$ and minimum weight $20$
was found in~\cite{GNW}.
Then $25$ more doubly even self-dual codes
 of length $120$ and minimum weight $20$
were found in~\cite{YW}.
From~\cite[Table~1]{YW},
these codes have different weight enumerators.
By a non-exhaustive search, we found $500$
doubly even self-dual four-circulant
codes of length $120$ and minimum weight $20$.
The numbers $a$ of codewords of weight $20$ in these codes
are listed in Table~\ref{Tab:WD120}.
It follows that 
these codes and the $26$ codes in~\cite{GNW} and \cite{YW}
have distinct weight enumerators.
Hence, we have the following:

\begin{prop}
There are at least $526$ inequivalent 
doubly even self-dual codes of length $120$ and minimum weight $20$.
\end{prop}

Our feeling is that the number of inequivalent 
doubly even self-dual codes  of length $120$ and minimum weight $20$
might be even bigger.

The first rows $r_A$ and $r_B$  of $A$ and $B$
in the generator matrices~\eqref{eq:GM}
of the $500$ codes are listed in 
\url{http://www.math.is.tohoku.ac.jp/~mharada/Paper/120-d20.txt}.
As an example, we list the first rows $r_A$ and $r_B$
for ten codes in Table~\ref{Tab:120}.

%%%%%%%%%%%%%%%%%%%%%%%%%%%%%%
\begin{table}[thbp]
\caption{Weight enumerators for length $120$}
\label{Tab:WD120}
\begin{center}
%{\small
%{\footnotesize
{\scriptsize
\begin{tabular}{l}
\noalign{\hrule height0.8pt}
\multicolumn{1}{c}{Numbers of codewords of weight $20$} \\ 
\hline
93180, 93936, 94512, 95136, 95202, 95376, 95496, 95532, 95826, 95946, 95952, 96012, 96096, 96126, \\
96156, 96216, 96240, 96312, 96336, 96360, 96366, 96372, 96486, 96540, 96576, 96666, 96690, 96720, \\
96762, 96780, 96816, 96840, 96846, 96876, 96906, 96912, 96936, 96996, 97026, 97056, 97092, 97116, \\
97176, 97230, 97260, 97266, 97272, 97296, 97326, 97356, 97422, 97446, 97452, 97476, 97566, 97572, \\
97590, 97596, 97626, 97632, 97656, 97716, 97746, 97770, 97776, 97782, 97836, 97842, 97866, 97890, \\
97896, 97926, 97950, 97962, 97986, 98016, 98040, 98076, 98130, 98136, 98166, 98196, 98220, 98226, \\
98250, 98256, 98262, 98286, 98292, 98316, 98346, 98412, 98466, 98496, 98502, 98526, 98532, 98556, \\
98562, 98580, 98586, 98610, 98616, 98622, 98640, 98646, 98670, 98676, 98682, 98700, 98706, 98712, \\
98730, 98742, 98772, 98796, 98802, 98826, 98832, 98856, 98886, 98910, 98916, 98940, 98952, 98976, \\
99000, 99036, 99066, 99090, 99096, 99120, 99126, 99156, 99162, 99180, 99186, 99210, 99216, 99222, \\
99240, 99246, 99252, 99270, 99282, 99306, 99312, 99330, 99336, 99342, 99372, 99390, 99396, 99402, \\
99432, 99450, 99456, 99486, 99516, 99540, 99546, 99576, 99612, 99666, 99672, 99690, 99696, 99702, \\
99720, 99726, 99750, 99756, 99786, 99792, 99810, 99816, 99846, 99876, 99906, 99936, 99942, 99966, \\
99972, 99996, 100026, 100032, 100062, 100086, 100110, 100116, 100122, 100140, 100146, 100170, \\
100176, 100182, 100200, 100206, 100212, 100236, 100242, 100260, 100266, 100290, 100296, 100350, \\
100356, 100380, 100446, 100452, 100476, 100482, 100500, 100506, 100512, 100536, 100542, 100560, \\
100566, 100590, 100596, 100626, 100650, 100656, 100662, 100680, 100686, 100716, 100722, 100746, \\
100752, 100770, 100776, 100782, 100800, 100806, 100842, 100860, 100872, 100896, 100902, 100920, \\
100926, 100956, 100980, 100986, 100992, 101046, 101052, 101070, 101076, 101082, 101106, 101112, \\
101130, 101136, 101142, 101160, 101166, 101196, 101202, 101226, 101232, 101250, 101256, 101280, \\
101286, 101316, 101376, 101382, 101400, 101406, 101412, 101436, 101442, 101472, 101496, 101526, \\
101532, 101550, 101556, 101586, 101616, 101622, 101640, 101646, 101652, 101670, 101676, 101700, \\
101706, 101730, 101736, 101760, 101766, 101772, 101790, 101796, 101802, 101820, 101826, 101850, \\
101856, 101862, 101880, 101892, 101910, 101916, 101940, 101946, 101952, 101970, 101976, 101982, \\
102000, 102006, 102030, 102036, 102042, 102066, 102072, 102096, 102120, 102126, 102150, 102156, \\
102180, 102186, 102210, 102216, 102240, 102246, 102252, 102270, 102312, 102336, 102342, 102360, \\
102366, 102372, 102402, 102420, 102426, 102456, 102480, 102486, 102492, 102516, 102540, 102546, \\
102570, 102576, 102582, 102606, 102636, 102660, 102666, 102672, 102690, 102696, 102702, 102726, \\
102732, 102750, 102756, 102780, 102786, 102792, 102816, 102840, 102846, 102870, 102876, 102906, \\
102930, 102936, 102942, 102966, 102972, 102996, 103002, 103020, 103026, 103032, 103050, 103056, \\
103080, 103086, 103092, 103116, 103140, 103146, 103176, 103182, 103206, 103236, 103266, 103272, \\
103296, 103320, 103326, 103332, 103356, 103380, 103386, 103410, 103416, 103422, 103452, 103500, \\
103506, 103530, 103560, 103566, 103590, 103596, 103632, 103650, 103656, 103686, 103692, 103710, \\
103716, 103722, 103740, 103746, 103752, 103770, 103776, 103800, 103806, 103830, 103836, 103860, \\
103896, 103932, 103962, 103986, 104022, 104046, 104076, 104106, 104166, 104220, 104226, 104232, \\
104256, 104286, 104316, 104346, 104436, 104442, 104496, 104502, 104532, 104556, 104580, 104592, \\
104616, 104622, 104646, 104652, 104676, 104736, 104772, 104796, 104820, 104880, 104886, 104892, \\
104910, 104916, 104970, 104982, 105066, 105096, 105156, 105336, 105396, 105426, 105456, 105510, \\
105546, 105576, 105636, 105666, 105696, 105762, 105966, 106152, 106236, 106266, 106290, 106386, \\
106626, 106662, 106812, 106836, 107220, 107406, 108486, 108600
\\
\noalign{\hrule height0.8pt}
\end{tabular}
}
\end{center}
\end{table}
%%%%%%%%%%%%%%%%%%%%%%%%%%%%%%%

%%%%%%%%%%%%%%%%%%%%%%%%%%%%%%
\begin{table}[thbp]
\caption{Doubly even self-dual codes  of length $120$ and minimum weight $20$}
\label{Tab:120}
\begin{center}
%{\small
{\footnotesize
%{\scriptsize
\begin{tabular}{c|c}
\noalign{\hrule height0.8pt}
$r_A$ & $r_B$ \\
\hline
$(100000111110000011111010001101)$&$(010000000010010010110101001111)$ \\
$(100001100101001111000100110010)$&$(100101100010011100001110011100)$ \\
$(100000011001010111001110101001)$&$(111101110101111111001100111010)$ \\
$(100001001111010001100000000011)$&$(010110111001100010000111011101)$ \\
$(100001011011001010010110000010)$&$(001010001100000100000010001000)$ \\
$(100000110010110111100100111110)$&$(100111000011110011101010001010)$ \\
$(100000000000111000011000111101)$&$(111000101010011000111000111011)$ \\
$(100001001111110101101101110111)$&$(011001000000100100101101110100)$ \\
$(100000001100001100111011101110)$&$(111111010001001110011000000000)$ \\
$(100000101111001000010100100110)$&$(001101011010111010101111011110)$ \\
\noalign{\hrule height0.8pt}

\end{tabular}
}
\end{center}
\end{table}
%%%%%%%%%%%%%%%%%%%%%%%%%%%%%%%

%%%%%%%%%%%%%%%%%%%%%%%%%%%%%%%%%%%%%%%%%%%%%%
\section{New doubly even self-dual codes
of length 128 and minimum weight 20}

%By Gleason's theorem (see~\cite{MS73}),
From \eqref{eq:WII},
the possible weight enumerator of a doubly even 
self-dual code  of length $128$ and minimum weight $20$
is determined as follows:
\begin{align*}
&
1 
+ a y^{20} 
+ (13228320  - 6 a) y^{24} 
+ (2940970496  - 89 a) y^{28} 
\\ &
+ (320411086380  + 1500 a) y^{32} 
+ (18072021808640  - 10925 a) y^{36} 
\\ &
+ (552523816524960  + 51186 a) y^{40} 
+ (9491115264030720  - 173451 a) y^{44} 
\\ &
+ (94116072808107840  + 449616 a) y^{48} 
+ (549827773219608576  - 920550 a) y^{52} 
\\ &
+ (1920594735166941760  + 1518100 a) y^{56} 
\\ &
+ (4051982995220321280  - 2040714 a) y^{60} 
\\ &
+ (5193576851944293670  + 2250664 a) y^{64}
+ \cdots, 
\end{align*}
where $a$ is the number of codewords of weight $20$.

The existence of a doubly even self-dual code
 of length $128$ and minimum weight $20$
is known~\cite{G}.
By a non-exhaustive search, we found $200$
doubly even self-dual four-circulant
codes  of length $128$ and minimum weight $20$.
These codes have distinct weight enumerators, where
the numbers $a$ of codewords of weight $20$ 
are listed in Table~\ref{Tab:WD128}.

\begin{prop}
There are at least $200$ inequivalent 
doubly even self-dual codes of length $128$ and minimum weight $20$.
\end{prop}

Our feeling is that the number of inequivalent 
doubly even self-dual codes  of length $128$ and minimum weight $20$
might be even bigger.  

The first rows $r_A$ and $r_B$  of $A$ and $B$
in the generator matrices~\eqref{eq:GM}
of the $200$ codes are listed in 
\url{http://www.math.is.tohoku.ac.jp/~mharada/Paper/128-d20.txt}.
As an example, we list the first rows $r_A$ and $r_B$
for ten codes in Table~\ref{Tab:128}.

%%%%%%%%%%%%%%%%%%%%%%%%%%%%%%
\begin{table}[thbp]
\caption{Weight enumerators for length $128$}
\label{Tab:WD128}
\begin{center}
%{\small
{\footnotesize
%{\scriptsize
\begin{tabular}{l}
\noalign{\hrule height0.8pt}
\multicolumn{1}{c}{Numbers of codewords of weight $20$} \\ 
\hline
21376, 21824, 22016, 22400, 22464, 22880, 22944, 23008, 23104, 23136, 23232, \\
23296, 23328, 23360, 23392, 23520, 23552, 23616, 23648, 23680, 23808, 23936, \\
24000, 24032, 24064, 24096, 24128, 24160, 24192, 24224, 24256, 24288, 24320, \\
24352, 24384, 24416, 24448, 24480, 24512, 24544, 24576, 24640, 24672, 24704, \\
24736, 24768, 24800, 24832, 24864, 24896, 24928, 24960, 24992, 25024, 25056, \\
25088, 25120, 25152, 25184, 25216, 25248, 25280, 25312, 25344, 25376, 25408, \\
25440, 25472, 25504, 25536, 25568, 25600, 25632, 25664, 25696, 25728, 25760, \\
25824, 25856, 25888, 25920, 25952, 25984, 26016, 26048, 26080, 26112, 26144, \\
26176, 26208, 26240, 26272, 26304, 26336, 26368, 26400, 26432, 26464, 26496, \\
26528, 26560, 26592, 26624, 26656, 26688, 26720, 26752, 26784, 26816, 26848, \\
26880, 26912, 26944, 26976, 27008, 27040, 27072, 27104, 27136, 27168, 27200, \\
27232, 27264, 27296, 27328, 27360, 27392, 27424, 27456, 27488, 27520, 27584, \\
27616, 27648, 27680, 27712, 27744, 27776, 27808, 27840, 27872, 27904, 27936, \\
27968, 28000, 28032, 28064, 28096, 28128, 28160, 28192, 28224, 28256, 28288, \\
28320, 28352, 28384, 28416, 28448, 28480, 28512, 28544, 28576, 28608, 28640, \\
28736, 28768, 28800, 28832, 28864, 28896, 28928, 28992, 29024, 29056, 29088, \\
29120, 29152, 29216, 29248, 29312, 29344, 29376, 29536, 29600, 29632, 29696, \\
29760, 29792, 29824, 29856, 29888, 30048, 30144, 30176, 30208, 30240, 30304, \\
30368, 31584\\
\noalign{\hrule height0.8pt}
\end{tabular}
}
\end{center}
\end{table}
%%%%%%%%%%%%%%%%%%%%%%%%%%%%%%%

%%%%%%%%%%%%%%%%%%%%%%%%%%%%%%
\begin{table}[thbp]
\caption{Doubly even self-dual codes  of length $128$ and minimum weight $20$}
\label{Tab:128}
\begin{center}
%{\small
{\footnotesize
%{\scriptsize
\begin{tabular}{c|c}
\noalign{\hrule height0.8pt}
$r_A$ & $r_B$ \\
\hline
$(10000111111110100010000110111100)$&$(00101010111100011101101011110100)$\\
$(10000111010101110101100011100110)$&$(10001000101110011011111011100110)$\\
$(10000010011000000001101010010001)$&$(01010001111001100001001011000010)$\\
$(10000111000111101111100111100111)$&$(11010000000111100110110101111111)$\\
$(10000010001100010100011010000111)$&$(01110100011000110110100110000110)$\\
$(10000111010101110010100010001010)$&$(01011110001110011011111011100111)$\\
$(10000100011111000010100010100010)$&$(10001101000000111010111010010101)$\\
$(10000101110011001010100011011010)$&$(00111111111010001111100110111000)$\\
$(10000000110010100000011111111000)$&$(11110000010011001110110110110101)$\\
$(10000010010000110011101110110100)$&$(10101000001000010101011011100001)$\\
\noalign{\hrule height0.8pt}
\end{tabular}
}
\end{center}
\end{table}
%%%%%%%%%%%%%%%%%%%%%%%%%%%%%%%

%%%%%%%%%%%%%%%%%%%%
\bigskip
\noindent
{\bf Acknowledgment.}
This work was supported by JSPS KAKENHI Grant Number 15H03633.
The author would like to thank
the anonymous referees for their useful comments on the manuscript.

%%%%%%%%%%%%%%%%%%%  References  %%%%%%%%%%%%%%%%%%%%%%%%

%%%%%%%%%%%%%%%%%%%%%%%%%%%%%%%%%


\begin{thebibliography}{99}

\bibitem{Magma}W. Bosma, J. Cannon and C. Playoust,
The Magma algebra system I: The user language,
{\sl J. Symbolic Comput.}
{\bf 24} (1997), 235--265.

\bibitem{BP} R.~Brualdi and V.~Pless,
Weight enumerators of self-dual codes,
{\sl IEEE Trans.\ Inform.\ Theory}
{\bf 37} (1991), 1222--1225.

\bibitem{C-S}J.H. Conway and N.J.A. Sloane,
{A new upper bound on the minimal distance of self-dual codes},
{\sl IEEE Trans.\ Inform.\ Theory}
{\bf 36} (1990), 1319--1333.

\bibitem{G} P. Gaborit, 
``Table of Type II codes'',
Online available at 
\url{http://www.unilim.fr/pages_perso/philippe.gaborit/SD/GF2/GF2II.htm},
Accessed on October 6, 2017.

\bibitem{GNW} P. Gaborit, C.-S. Nedeloaia and W. Wassermann, 
On the weight enumerators of duadic and quadratic residue codes,
{\sl IEEE Trans.\ Inform.\ Theory}
{\bf 51} (2005),  402--407.

\bibitem{Gleason}A.M. Gleason, 
Weight polynomials of self-dual codes and the MacWilliams identities,
Actes du Congr\'es International des Math\'ematiciens (Nice, 1970), Tome 3, 
pp.~211--215. Gauthier-Villars, Paris, 1971. 

\bibitem{H112} M. Harada,
An extremal doubly even self-dual code of length $112$,
{\sl Electron.\ J. Combin.}
{\bf 15} (2008), Note 33, 5 pp. 

\bibitem{4cir}
M. Harada, W. Holzmann, H. Kharaghani and M. Khorvash, 
Extremal ternary self-dual codes constructed from negacirculant matrices,
{\sl Graphs Combin.}
{\bf  23} (2007), 401--417.

% \bibitem{HM-CR} M. Harada and A. Munemasa,
% On the covering radii of extremal doubly even self-dual codes,
% {\sl Adv.\ Math.\ Commun.}
% {\bf 1} (2007), 251--256.

\bibitem{MS73} C.L. Mallows and N.J.A. Sloane,
{An upper bound for self-dual codes},
{\sl Inform.\ Control}
{\bf 22} (1973), 188--200.

% \bibitem{Rains} E.M.~Rains,
% {Shadow bounds for self-dual codes},
% {\sl IEEE Trans.\ Inform.\ Theory}
% {\bf 44} (1998), 134--139.

\bibitem{RS-Handbook} E.\ Rains and N.J.A.\ Sloane,
{``Self-dual codes,''} {Handbook of Coding Theory},
V.S. Pless and W.C. Huffman (Editors),
Elsevier, Amsterdam 1998, pp.\ 177--294.

\bibitem{YW} R. Yorgova and A. Wassermann, 
Binary self-dual codes with automorphisms of order $23$,
{\sl Des.\ Codes Cryptogr.}
{\bf 48} (2008), 155--164.

\end{thebibliography}
\end{document}